\documentclass{amsart}
\usepackage{amsmath}
\usepackage{amsfonts}
\usepackage{amssymb}
\usepackage{mathtools}
\usepackage{thmtools}
\usepackage{thm-restate}
\usepackage{tikz}
\usepackage{multirow}
\usepackage{url}
\usepackage{array}
\usepackage{float}
\usepackage{tikz-cd}

\usepackage[unicode,pdfborder={0 0 0}]{hyperref}
\usepackage[english]{cleveref}
\usepackage{amsrefs}
\usepackage{mathrsfs}

\newtheorem{prop}{Proposition}[section]
\newtheorem{lemma}[prop]{Lemma}
\newtheorem{teo}[prop]{Theorem}

\newtheorem{cor}[prop]{Corollary}
\newtheorem{defi}[prop]{Definition}

\newtheorem{rmk}[prop]{Remark}

\newcommand{\IP}{\mathbb{P}}
\newcommand{\IC }{\mathbb{C}}

\newcommand{\IZ}{\mathbb{Z}}

\newcommand{\rk}{\mathrm{rk}}

\newcommand{\deco}{\IP(A)\cap G(3,V_6)=\emptyset}

\DeclareMathOperator{\Aut}{Aut}

\DeclareMathOperator{\Bir}{Bir}

\DeclareMathOperator{\coker}{Coker}

\DeclareMathOperator{\dete}{det}

\DeclareMathOperator{\GL}{GL}

\DeclareMathOperator{\id}{id}

\DeclareMathOperator{\NS}{NS}
\DeclareMathOperator{\ns}{NS}
\DeclareMathOperator{\PGL}{PGL}

\DeclareMathOperator{\PSL}{PSL}

\DeclareMathOperator{\sing}{Sing}

\DeclareMathOperator{\SO}{SO}
\DeclareMathOperator{\Sp}{Sp}

\DeclareMathOperator{\Spin}{Spin}

\DeclareMathOperator{\Sym}{Sym}

\DeclareMathOperator{\LG}{\mathbb{LG}(\bigwedge^3 V_6)}
\newcommand*{\sheafhom}{\mathcal{H}\kern -.5pt om}

\tikzset{
    labl/.style={anchor=south, rotate=90, inner sep=.5mm}
}

\makeatletter
\def\blfootnote{\xdef\@thefnmark{}\@footnotetext}

\begin{document}

\title{Double EPW sextics associated to Gushel-Mukai surfaces}

\author{Pietro Beri}
\address{\parbox{0.9\textwidth}{
Universit\'e Paris Cit\'e and Sorbonne Universit\'e, CNRS, IMJ-PRG, F-75013 Paris, France
\vspace{1mm}}}
\email{{beri@imj-prg.fr }}

\thanks{The author is supported by the ERC Synergy Grant ERC-2020-SyG-854361-HyperK}

\begin{abstract}
Works by O'Grady allow to associate to a 2-dimensional Gushel-Mukai variety, which is a K3 surface, a double $EPW$ sextic. 
We characterize the $K3$ surfaces whose associated double $EPW$ sextic is smooth.
As a consequence, we are able to produce symplectic actions on some families of smooth double $EPW$ sextics which are hyper-K\"ahler manifolds.

We also provide bounds for the automorphism group of Gushel-Mukai varieties in dimension 2 and higher.
\end{abstract}

\maketitle

\section{Introduction}

Double $EPW$ sextics are an important family of hyper-K\"ahler manifolds i.e. compact simply connected
K\"ahler manifolds with a unique, up to scalar, holomorphic two-form, which is everywhere non-degenerate.
The linebase example of hyper-K\"ahler manifolds are $K3$ surfaces; a classical example are double covers of $\IP^2$.
Double $EPW$ sextics are a generalization of it, as they come with a structure of double covers of special sextic hypersurfaces in $\IP^5$, the so-called $EPW$ sextics.

\medskip

In Section \ref{le EPW e tut quant} we present the main objects that come into play.
The basis of the theory of double $EPW$ sextics has been developed by O'Grady in an influential series of papers.
Already from \cite{ogrady2006a}, an important connection is observed between (double) $EPW$ sextics and ordinary $2$-dimensional Gushel-Mukai ($GM$) varieties, complete intersections of a linear space and a quadric hypersurface inside a Grassmannian.
This link has been successively extended to higher-dimensional $GM$ varieties by Iliev and Manivel in \cite{iliev2011}
and then developed in detail by Debarre and Kuznetsov, in a series of papers that lay the groundwork for an extensive study of the beautiful and intricate interplay between the two families.

\medskip

A natural question about this relation is: \emph{can we give conditions on a $GM$ variety to be associated to a double $EPW$ sextic which is smooth?}
Indeed, whenever the double cover is smooth, it is a hyper-K\"ahler manifold.

\medskip

In Section \ref{mio divisoriale} we answer this question for $2$-dimensional $GM$ varieties, which are Brill-Noether general $\langle 10 \rangle$-polarized $K3$ surfaces. We provide an answer in terms of geometry and in terms of period: the following result summarizes \Cref{mio risultato su dominio: cor} and \Cref{mio risultato su dominio: versione pic}.
\begin{teo}
\label{teo di intro}
The double $EPW$ sextic associated to a strongly smooth $K3$ surface $S=\IP^6\cap G(2,V_5)\cap Q$ is smooth if and only if $S$ contains neither lines nor quintic elliptic pencils.
Equivalently, a $\langle 10 \rangle$-polarized $K3$ surface is associated to a double $EPW$ sextic which is a hyper-K\"ahler manifold if and only if it does not lie in six (explicitly described) divisors in the corresponding moduli space.
\end{teo}

This result relies on the tools provided by Debarre and Kuznetsov which we presented in Section \ref{le EPW e tut quant} as well as on a careful description of the involution acting on the moduli space of $EPW$ sextics, described for the first time by O'Grady \cite{ogrady2006}.
It is interesting to note that smoothness was already known in the very general case, but \Cref{teo di intro} holds without any general assumption: eliminating the hypothesis of very generality is usually a challenging problem and we believe that our result is interesting in this spirit.  

\medskip

As an application of \Cref{teo di intro}, in Section \ref{auto indotti} we produce symplectic actions for various groups on families of double $EPW$ sextics which are hyper-K\"ahler,
from lattice-theoretic considerations on automorphisms of $K3$ surfaces.

\medskip

Section \ref{automorfismi vari e eventuali} is devoted to finding bounds for the automorphism group 
of $GM$ varieties and on actions of groups on $EPW$ sextics. For example, in \Cref{gli auto di una GM che e una K3} we show that the automorphism group of $S$ is a finite subgroup of $PGL(2,\IC)$ and can only be symplectic; we also provide some results for $GM$ varieties in higher dimension. Most of these results are obtained by studing the automorphism group of Fano varieties containing a $GM$ variety.

\medskip

A previous version of this paper contained some results about lifting automorphisms of $EPW$ sextics to their double cover; in the meanwhile, a much stronger result has been obtained by Kuznetsov in \cite{debarre2021}, see \Cref{tuuuto symply}. We used this result to simplify the proof of \Cref{famiglie defo di K3 GM simply} and to prove the symplecticity of the automorphism group of a strongly smooth $K3$ surface in $\IP^6$, see \Cref{gli auto di una GM che e una K3}.

\medskip

\textbf{Acknowledgments:} Most of the results of this paper are part of my PhD thesis. I am grateful to Alessandra Sarti for her constant support, for the supervision of this work and for many useful discussions.
I would like to thank Giovanni Mongardi, Enrico Fatighenti and Simone Novario for many interesting discussions. Finally I would also like to thank the anonymous referee for many helpful comments and suggestions.

\subsection{Notations}

Let $X$ be a topological space. A property holds for $x\in X$ general if the condition is satisfied by all the points inside an open subset of $X$.
A property holds for $x\in X$ \emph{very} general if the condition is satisfied by all the points in the complement of a countable union of closed subspaces inside $X$.

\medskip
Given a complex vector space $V$, the Grassmannian of $k$-dimensional vector subspaces in $V$ will be denoted by $G(k,V)$; every Grassmannian we consider will be embedded in the projective space $\IP(\bigwedge^k V)$ via the Pl\"ucker embedding. Non-zero decomposable vectors in $\IP(\bigwedge^k V)$ are elements that lie in $G(k,V)$.

\begin{defi}
\label{quadriche pfaffiane}
Consider $V_5\cong\IC^5$. Following \cite{iliev2011}, for $v\in V_5-\lbrace 0 \rbrace$ we call \emph{Pfaffian quadric} $P_v$ the symmetric bilinear form $x \mapsto vol (v\wedge x \wedge x )$ on $\bigwedge^2 V_5$.
\end{defi}
The Grassmannian is the intersection of all the Pfaffian quadrics;
since $\ker(P_v)=v\wedge V_5$, every Pfaffian quadric $P_v$ has rank $6$.

An automorphism $\alpha$ on a hyper-K\"ahler manifold $X$ is \emph{symplectic} if it acts trivially on a $2$-form on $X$, \emph{non-symplectic} otherwise. We always consider $H^2(X,\IZ)$ endowed with the Beauville-Bogomolov (BB) form. A $\langle 2t \rangle$-polarization is a polarization whose square is $2t$ with respect to the BB form.\\
We write $S^{[2]}$ for the Hilbert square on a $K3$ surface $S$. For any $H\in \ns(S)$, we denote by $H_2$ the induced class on $S^{[2]}$ and by $\delta\in \ns(S^{[2]})$ the class of the divisor such that $2\delta$ parametrizes non-reduced subschemes.

\section{EPW sextics and Gushel-Mukai varieties}
\label{le EPW e tut quant}

\subsection{$EPW$ sextics}
\label{inizio di epw sestiche}

O'Grady provided slightly different points of view on the construction of Eisenbud-Popescu-Walter sextics, see for example \cite{ogrady2006a}, \cite{ogrady2010} and \cite{ogrady2012}; here we follow mainly \cite{ogrady2012}, with a view to \cite{debarre2019a}, by Debarre and Kuznetsov.

Let $V_6$ be a $6$-dimensional vector space over $\mathbb{C}$, on which we fix a volume form $\bigwedge^6 V_6\xrightarrow{\sim} \IC$, which in turns induces a symplectic form $\omega$ on $\bigwedge^3 V_6$.

\medskip

In \cite[Example 9.3]{eisenbud2001}, Eisenbud, Popescu and Walter introduced the Lagrangian subbundle $F$ of $\mathcal{O}_{\IP(V_6)}\otimes \bigwedge^3 V_6$, whose fiber over $[v]$ is $F_v=v\wedge \bigwedge^2 V_6$.

Given $v\in V_6-\lbrace 0 \rbrace$, we can fix a decomposition $V_6\cong \mathbb{C}v\oplus V_5$ for some hyperplane $V_5\subset V_6$, which induces $\bigwedge^3 V_6\cong \bigwedge^3 V_5 \oplus F_v$. Every element of $F_v$ can be written in the form $v\wedge \eta$ for some $\eta\in \bigwedge^2 V_5$ and this induces an isomorphism of vector spaces
\begin{align}
\label{decomposizione V6}
\rho\colon F_v  &\xrightarrow{\sim} \bigwedge^2 V_5\\
 v\wedge \eta &\mapsto   \eta.
\end{align}
This observation will be useful later.

\medskip

We call $\mathbb{LG}(\bigwedge^3 V_6)\subset G(10,\bigwedge^3 V_6)$ the \emph{symplectic Grassmannian}, parametrizing Lagrangian subspaces with respect to the symplectic form on $\bigwedge^3 V_6$; since two volume forms differ by a non-zero constant, $\mathbb{LG}(\bigwedge^3 V_6)$ does not depend on the choice of the volume form.
From now on, $A\in \mathbb{LG}(\bigwedge^3 V_6)$ will be a Lagrangian subspace in $\bigwedge^3 V_6$ such that $\deco$. 

\medskip

For a fixed $A$, we consider the map of vector bundles $\lambda_A\colon F\to \mathcal{O}_{\mathbb{P}(V_6)} \otimes A^\vee$
such that $(\lambda_A)_{[v]}(x)=\omega(x,-)_{|_{A}}$ for every $[v]\in \IP(V_6)$.
Since $A$ is a Lagrangian subspace, the map is injective on $[v]\in \IP(V_6)$ if and only if $F_v\cap A=0$.
Since $\rk(F)=\rk(\mathcal{O}_{\mathbb{P}(V_6)}\otimes A^\vee)$, we can consider the map
$ \dete (\lambda_A)\colon \dete (F)\to \dete  (\mathcal{O}_{\mathbb{P}(V_6)}\otimes A^\vee).$ 
\begin{defi}
\label{definizione di epw sestica}
The determinantal variety $Y_A= Z( \dete (\lambda_A))$ is a sextic hypersurface (see \cite[(1.8)]{ogrady2006a}), called \emph{Eisenbud Popescu Walter ($EPW$) sextic}.
\end{defi}

We also associate to $A$ a stratification of $\mathbb{P}(V_6)$:
for $k\geq 0$ we define
\begin{equation}
\label{definizione delle varie epw varieta}
Y^{\geq k}_A=\lbrace [v] \in \mathbb{P}(V_6) | \dim(F_v\cap A)\geq k \rbrace;
\end{equation}
\[Y^{k}_A=\lbrace [v] \in \mathbb{P}(V_6) | \dim(F_v\cap A)= k \rbrace.\]
The sets $Y^{\geq k}_A$ are degeneracy loci and are thus endowed with a natural structure of scheme.
For every $k$, the variety $Y^{\geq k+1}_A$ is a closed subvariety of $Y^{\geq k}_A$ and clearly $Y_A=Y_A^{\geq 1}$.
\begin{prop}
\label{mult 4 e decomponibili}
\label{sing di Y_A}
\label{sing della superficie}
\label{mult 3 se non ho decomponibili}
\label{superficie Y^2_A}
If $\deco$, then $Y_A$ is integral and normal and
$\sing(Y_A)=Y^{\geq 2}_A$.
Moreover $Y^{\geq 2}_A$ is a normal integral surface whose singular locus is $Y^{3}_A$, which in turn is finite; when smooth, the surface $Y^{\geq 2}_A$ is of general type.
Finally $Y^{\geq 4}_A=\emptyset$.
\begin{proof}
Integrality and normality for $Y_A$ follow from $\sing(Y_A)=Y^{\geq 2}_A$.
For the rest, see \cite[Corollary $2.5$]{ogrady2012}, \cite[Proposition $2.9$]{ogrady2012}, \cite[Claim $3.7$]{ogrady2010}, \cite[Theorem $B.2$]{debarre2015} and \cite[Proposition $1.10$]{ferretti2012}.
\end{proof}
\end{prop}

The situation can be much more complicated for $EPW$ sextics associated to a Lagrangian subspace $A$ which contains a non-zero decomposable vector, however for our scope we will not need to deal with them.

\subsection{The dual EPW sextic}
One interesting feature of $EPW$ sextics is that they admit a dual counterpart.
The volume form on $V_6$ induces a volume form on $ V^\vee_6$ too, and thus a symplectic form on $\bigwedge^3 V^\vee_6\cong (\bigwedge^3 V_6)^\vee$, hence
for $A\in \mathbb{LG}(\bigwedge^3 V_6)$ we also have
\[A^\perp=\lbrace \phi\in (\bigwedge^3 V_6)^\vee|\mbox{ } \phi(a)=0 \mbox{ for every }a\in A\rbrace\in \mathbb{LG}((\bigwedge^3 V_6)^\vee);\]
$\deco$ if and only if $\IP(A^\perp)\cap G(3,V^\vee_6)=\emptyset$, see \cite[Section $2.6$]{ogrady2012}.
\begin{defi}
We call \emph{dual $EPW$ sextic} the $EPW$ sextic associated to $A^\perp$.
\end{defi}

When $\deco$, the hypersurfaces $Y_A\subset \mathbb{P}(V_6)$ and $Y_{A^\perp}\subset \mathbb{P}(V^\vee_6)$ are projectively dual;
this has been proved by O'Grady, see also \cite[Proposition $B.3$]{debarre2015}.

\medskip

As showed in \cite[Section $2.6$]{ogrady2012}, the $EPW$ strata associated to $A^\perp$ can also be described as
\[Y^{\geq k}_{A^\perp}=\lbrace [V_5] \in \mathbb{P}(V^\vee_6) | \dim(\bigwedge^3 V_5 \cap A)\geq k   \rbrace\]
\begin{defi}
\label{i divisori sulle lagrangiane}
\label{ultimo divisore nello spazio della lagrangiane}
\label{divisori in LG e LG duale stesso spazio per PI}
We denote by 
\[\Sigma=\lbrace A\in \mathbb{LG}(\bigwedge^3 V_6)\mbox{ } |\mbox{ } \mathbb{P}(A)\cap G(3,V_6)\neq \emptyset \rbrace \]
the set of Lagrangian subspaces $A$ admitting a non-zero decomposable vector, by
\[\Delta=\lbrace A\in \mathbb{LG}(\bigwedge^3 V_6)\mbox{ } | \mbox{ } Y^{\geq 3}_A\neq \emptyset\rbrace\] 
the set of Lagrangian subspaces whose third stratum is not empty, and by 
\[\Pi=\lbrace A\in \mathbb{LG}(\bigwedge^3 V_6)\mbox{ } | \mbox{ } Y^{\geq 3}_{A^\perp}\neq \emptyset\rbrace \]
the set of Lagrangian subspaces $A\in \LG$ such that the third stratum associated to the dual Lagrangian $A^\perp$ is not empty.\\
Following O'Grady, we write $\mathbb{LG}(\bigwedge^3 V_6)^0=\mathbb{LG}(\bigwedge^3 V_6)- (\Delta\cup \Sigma)$.

\end{defi}
The subsets $\Sigma$, $\Delta, \Pi$ descend to distinct irreducible divisors in $\mathbb{LG}(\bigwedge^3 V_6)//PGL(6)$, see \cite[Proposition $3.1$ Item $1)$]{ogrady2012} and \cite[Proposition $2.2$]{ogrady2010}. We denote their quotients in the same way.

\medskip

\begin{prop}\cite[Proposition $B.8$]{debarre2015}
\label{auto fissano A}
If $A\notin \Sigma$, the automorphism group of $Y_A$ is $\lbrace \alpha\in \PGL(V_6)$ such that $ \left(\bigwedge^3 \alpha\right)(A)=A \rbrace$. In particular every $\alpha\in \Aut(Y_A)$ fixes $Y_A^k$ for every $k\geq 0$.
\end{prop}

From now on we always identify $PGL(V_6)$ and $PGL(V^\vee_6)$ through the natural isomorphism sending
$\phi\in \PGL(V_6)$ to $(\phi^{-1})^\vee$, which maps $[V_5]$ to $[\phi(V_5)]$. 

\begin{cor}
\label{auto di EPW e suo duale sono uguali}
If $A\notin \Sigma$, then
$\Aut(Y_A)=\Aut(Y_{A^\perp})$.
\end{cor}

\subsection{Double $EPW$ sextics}
\label{double covers per lagrangiani}

The importance of $EPW$ sextics stems mostly from the fact that, in the general case, they admit a ramified double cover which is a hyper-K\"ahler manifold, as was proved by O'Grady in \cite{ogrady2006a}.
We also refer to \cite[Theorem $5.2$, Item $1)$]{debarre2019a}.

We denote by $\mathcal{R}=\coker (\lambda_A)_{|_{Y_A}}$ the \emph{first Lagrangian cointersection sheaf}; it is a reflexive sheaf, see \cite[Proposition $4.3$, Item $1)$]{ogrady2006a}.
\begin{teo}(O'Grady)
\label{esistenza della double EPW sestica}
Consider $A\in \mathbb{LG}(\bigwedge^3 V_6)$ such that $\deco$. Then there is a unique double cover $f_A\colon X_A\to Y_A$ with branch locus $Y^{\geq 2}_A$ such that
\begin{equation}\label{fasci per double EPW}
(f_A)_*\mathcal{O}_{X_A}\cong \mathcal{O}_{Y_A} \oplus \mathcal{R}(-3),
\end{equation}
The variety $X_A$ is integral and normal, and its singular locus is $f^{-1}_A(Y^3_A)$.
\end{teo}
A double cover of $Y_A$ can be constructed even when $\IP(A)\cap G(3,V_6)\neq \emptyset$, provided that $Y_A\neq \IP(V_6)$. However, in this case $X_A$ is never smooth.
\begin{defi}
We call \emph{double $EPW$ sextic} the double cover $X_A$, and we denote by $\iota_A$ the associated covering involution.
\end{defi}

\begin{teo}\cite[Theorem $4.25$]{ogrady2010}
\label{double EPW sextic hyperkahler}
Suppose that $X_A$ is smooth i.e. $A\in \LG^0$. Then $X_A$ is a hyper-K\"ahler manifold equivalent by deformation to the Hilbert square on a K3 surface. The involution $\iota_A$ is non-symplectic. 
\end{teo}

The ample class $D_A=f^*_A\mathcal{O}_{Y_A}(1)\in \ns(X_A)$
is the only primitive polarization on $X_A$ coming from $Y_A$, by Lefschetz hyperplane theorem.
It has square $2$ with respect to the BB form;
the general element, inside the moduli space $\mathcal{M}_2$ (see \cref{subsec involuzione} for a definition)
is isomorphic to $(X_A,D_A)$ for some $A\in \LG^0$. 

\medskip

We denote by $\Aut_{D_A}(X_A)\subseteq \Aut(X_A)$ the group of automorphisms fixing $D_A$, or equivalently commuting with $\iota_A$.
There is a short exact sequence
\begin{equation}\label{sequenza lift alla doubleEPW}
1\to \lbrace \id,\iota_A \rbrace\to \Aut_{D_A}(X_A) \to \Aut(Y_A) \to 1.
\end{equation}

\begin{prop}\cite[Proposition A.2]{debarre2021}
\label{tuuuto symply}
The sequence \eqref{sequenza lift alla doubleEPW} splits, so $\Aut_{D_A}(X_A)\cong \Aut(Y_A)\times \lbrace \id,\iota_A \rbrace$. 
When $X_A$ is smooth, under this isomorphism $\Aut(Y_A)$ is the group of symplectic automorphisms fixing the polarization $D_A$. 
\end{prop}

\subsection{Gushel-Mukai varieties}
\label{intro GM}

In this section we introduce Gushel-Mukai varieties and in the next one we explain how they admit an associated $EPW$ sextic.\\
We fix a $5$-dimensional complex vector space $V_5$.
For the following definition and results we refer to \cite[Definition $2.1$]{debarre2015}.
\begin{defi}
\label{defi GM}
Consider a vector subspace $W\subseteq \bigwedge^2 V_5$ of dimension $\dim W\geq 6$, and a quadric hypersurface $Q\subset \mathbb{P}(W)$. We call \emph{ordinary Gushel-Mukai ($GM$) intersection} the scheme
\[ Z= \mathbb{P}(W)\cap G(2,V_5)\cap Q . \]
it is an \emph{ordinary $GM$ variety} if $Z$ is integral 
and $\dim(Z)=\dim(W)-5$.
\end{defi}
When $Z$ is a $GM$ variety, it is a complete intersection inside $G(2,V_5)$ and has degree $10$ in $\IP(W)$. 

If $\dim(Z)\geq 3$, $Z$ is a Fano variety of index $n-2$ (see for example \cite[Theorem $2.3$]{debarre2015}). If $\dim(Z)=2$, let $H$ be the polarization given by $\mathcal{O}_{\IP(W)}(1)$: the polarized variety $(Z,H)$ is a Brill-Noether general $K3$ surface, see \cite{mukai}, also \cite[Theorem $10.3$]{knutsen2004}. The converse holds under a technical condition, which is strongly smoothness.

\begin{defi}\label{gushel mukai strongly smooth}
Let $Z=\IP(W)\cap G(2,V_5)\cap Q$ be an ordinary $n$-dimensional $GM$ variety. The \emph{Grassmannian hull of $Z$} is the $(n+1)$-dimensional intersection $M_Z=\IP(W)\cap G(2,V_5)$.
We say that $Z$ is \emph{strongly smooth} if both $Z$ and $M_Z$ 
are dimensionally transverse and smooth.
\end{defi}

Although $GM$ curves exist, here we never deal with them, since they are never strongly smooth or, equivalently, their associated double $EPW$ sextics are never smooth, as we will see later (\Cref{GM e EPW data}).

If $Z$ is strongly smooth, its Grassmannian hull has Picard rank 1 by the Lefschetz hyperplane theorem.
For every Brill-Noether general $K3$ surface $(S,H)$, the projective model $\phi_{|H|}(S)$ is a smooth ordinary $GM$ variety of dimension $2$, provided that it is strongly smooth (see \Cref{SS GM of dim 2 are ordinary} below), whereas smoothness and strongly smoothness are equivalent when $n\geq 3$.

\begin{defi}
\label{defi auto GM variety}
An isomorphism between $Z=\mathbb{P}(W)\cap G(2,V_5)\cap Q$ and $Z'=\mathbb{P}(W')\cap G(2,V'_5)\cap Q'$ is a linear map $\phi \colon \mathbb{P}(W)\to \mathbb{P}(W')$ such that
$\phi(Z)=Z'$.
We denote by $\Aut(Z,\IP(W))$ the group of automorphisms of $Z$.
\end{defi}

By the Lefschetz hyperplane theorem, a smooth $GM$ variety $Z$ of dimension at least $3$ has Picard rank $1$. So, in this case, $\Aut(Z,\IP(W))$ is the whole automorphism group of $Z$ as an abstract variety.

\medskip

Debarre and Kuznetsov provided an intrinsic characterization of normal $GM$ varieties in \cite[Theorem $2.3$]{debarre2015}; this leads to the definition of $GM$ data, which are a set of linear data which can be associated to any normal $GM$ variety. These collections of  objects are very useful to handle. 

\begin{defi}

\emph{Ordinary $GM$ data} $(W,V_6,V_5,\mu,q,\epsilon)$ of dimension $n$ consists of (we set $L=(V_6/V_5)^\vee$ for readability)
\label{definizione dei GM fuckin data}
\begin{itemize}
\item a $(n+5)$-dimensional vector space $W$;
\item a $6$-dimensional vector space $V_6$;
\item a hyperplane $V_5$ of $V_6$;
\item an \emph{injective} linear map $\mu\colon  W\otimes L \hookrightarrow \bigwedge^2 V_5$;
\item a linear map $q\colon V_6\to S^2 W^\vee$;
\item a linear isomorphism $\epsilon\colon \bigwedge^5 V_5\to L^{\otimes 2}$,
\end{itemize}
such that, for all $v\in V_5$, $w_1, w_2\in W$,
\begin{equation}
\label{cosa sensata GM data}
q(v)(w_1,w_2)=\epsilon(v\wedge \mu(w_1)\wedge \mu(w_2)).
\end{equation}
\end{defi}

\begin{defi}
An isomoprhism between two ordinary $GM$ data sets $(W,V_6,V_5,\mu,q,\epsilon)$ and $(W',V'_6,V'_5,\mu',q',\epsilon')$ is a triple of linear isomorphisms $\phi_W\colon W\to W'$, $\phi_V\colon V_6\to V'_6$, $\phi_L\colon L\to L'$ such that $\phi_V(V_5)=V'_5$, $\epsilon'\circ \bigwedge^5 \phi_V=\phi^{\otimes 2}_L\circ \epsilon$, and the following diagrams commute
\begin{equation}\label{diagrammi auto GM data}
\begin{tikzcd}
V_6  \ar[d,swap,"\phi_V"]\ar[rr,"q"] & & \Sym^2 {W}^\vee    \\
V'_6 \ar[rr,"q' "]  & & \Sym^2 {(W')}^\vee   \ar[u,swap,"\Sym^2 \phi^\vee_W"]
\end{tikzcd}
\quad \quad \quad \quad
\begin{tikzcd} 
W\otimes L\ar[r,"\mu"]\ar[d,swap,"\phi_W\otimes \phi_L"] & \bigwedge^2 V_5 \ar[d,"\bigwedge^2 \left({\phi_V}_{|_{V_5}}\right)"]  \\
W'\otimes L' \ar[r,"\mu'"]& \bigwedge^2 V'_5.
\end{tikzcd}
\end{equation}
\end{defi}

Debarre and Kuznetsov showed in \cite[Section $2.1$]{debarre2015} how to associate a set of $GM$ data $(W,V_6,V_5,\mu,q,\epsilon)$ to a normal $GM$ variety and to $GM$ data a $GM$ intersection which a priori may not be a normal variety; when it is the case, the two constructions are mutually inverse and behave well with respect to the definition of isomorphisms. To keep the exposition readable, we refer to their paper for the explicit correspondence.

\medskip

Finally, we point out that another class of $GM$ intersections, \emph{special $GM$ intersections}, exist, see for example \cite[Section $2.5$]{debarre2015}.

\subsection{A correspondence between data sets}
\label{corrispondenze}

Debarre and Kuznetsov defined another set of data. Everytime we consider a $6$-dimensional complex vector space $V_6$, we endow it with a volume form, which induces a symplectic form on $\bigwedge^3 V_6$. 

\begin{defi}
An \emph{ordinary Lagrangian data} is a collection $(V_6,V_5,A)$, where
\begin{itemize}
\item $V_6$ is a $6$-dimensional complex vector space,
\item $V_5\subset V_6$ is a hyperplane,
\item $A\in \LG$ is a Lagrangian subspace.
\end{itemize} 
The ordinary Lagrangian data $(V_6,V_5,A)$ and $(V'_6,V'_5,A')$ are isomorphic if there is a linear isomorphism $\phi:V_6\to V'_6$ such that $\phi(V_5)=V'_5$ and $(\bigwedge^3 \phi)(A)=A'$.
\end{defi}

The next result, by Debarre and Kuznetsov, sums up the results by O'Grady, \cite[Section $5$ and Section $6$]{ogrady2006a}, \cite[Section $4$]{ogrady2010}, in dimension $2$ and Iliev and Manivel, \cite[Section $2$]{iliev2011} in dimension $5$ and \cite[Section $4$]{iliev2011} in dimension $3,4$.

\begin{teo}
\label{GM e EPW data}
\label{strongly smooth e no decomponibili}
\label{strongly smooth e no decomponibili: cor}
For $n\in \lbrace 1,\ldots,5 \rbrace$,
there is a bijection between the set of isomorphism classes of Lagrangian data sets $(V_6,V_5,A)$, with
\begin{itemize}
\item $\deco$
\item $[V_5]\in Y^{5-n}_{A^\perp}$
\end{itemize}
and isomorphism classes of strongly smooth ordinary $GM$ \emph{varieties} of dimension $n$.
In particular, there are no strongly smooth ordinary $GM$ curves.

\begin{proof}
See \cite[Theorem $3.6$]{debarre2015} and \cite[Theorem $3.16$]{debarre2015}. For $n=1$, recall that $A^\perp$ admits non-zero decomposable elements if and only if $A$ does. Then we only need to observe that $Y^{\geq 4}_{A^\perp}=\emptyset$ by \Cref{mult 4 e decomponibili}.
\end{proof}
\end{teo}

It is interesting to note that the construction works for non strongly smooth $GM$ varieties, but in that case $A$ always admits non-zero decomposable elements.

\begin{defi}
Given some Lagrangian data $(V_6,V_5,A)$, the \emph{associated ordinary $GM$ variety} is the $GM$ intersection $Z$ obtained by \Cref{GM e EPW data} from $(V_6,V_5,A)$.
When not otherwise specified, we write $Z=\IP(W)\cap G(2,V_5)\cap Q$ where $Q$ is the projective quadric associated to $q(x)$ for any $x\in V_6-V_5$.

Given a strongly smooth ordinary $GM$ variety $Z$
we denote by $A(Z)$ the Lagrangian subspace associated to $Z$ by \Cref{GM e EPW data}.
\end{defi}

\begin{rmk}\label{SS GM of dim 2 are ordinary}
\normalfont
Strongly smooth special $GM$ varieties of dimension $2$ do not exist by \cite[Remark $3.17$]{debarre2015}, and \Cref{mult 4 e decomponibili}.
\end{rmk}

For any $(V_6,V_5,A)$ and $p\geq 2$, there is a short exact sequence
\begin{equation}\label{sequenza spazi}
0 \to \bigwedge^p V_5 \to \bigwedge^p V_6 \xrightarrow{\lambda_p}  \bigwedge^{p-1} V_5\otimes L^\vee  \to 0
\end{equation}

The following two remarks will be useful in the proof of \Cref{mio risultato su dominio: cor} in the frame of the study of 2-dimensional GM varieties:
although the relation between the latter and $EPW$ sextics is already clear in \cite{ogrady2006a} and \cite{ogrady2010}, the use of $GM$ data allows a very precise study.

\begin{rmk}
\label{primo remark}
\normalfont
Given some Lagrangian data $(V_6,V_5,A)$, \Cref{GM e EPW data} allows to describe explicitly the associated ordinary $GM$ variety $Z\subseteq \IP(\bigwedge^2 V_5)$. We fix $x\in V_6-V_5$, hence an isomorphism 
$\bigwedge^3 V_6 \cong \bigwedge^3 V_5\oplus F_x$.
Under this isomorphism and \eqref{decomposizione V6},
the sequence \eqref{sequenza spazi} for $p=3$ becomes
\begin{equation}\label{sequenza spazi vettoriali}
 0 \to  \bigwedge^3 V_5 \to  \bigwedge^3 V_6\cong \bigwedge^3 V_5\oplus F_x  \xrightarrow{\lambda_3} \bigwedge^{2} V_5 \to 0 .
\end{equation}
For $p=4$ the sequence \eqref{sequenza spazi} can be rewritten as $ 0 \to  \bigwedge^4 V_5 \to  \bigwedge^4 V_6  \xrightarrow{\lambda_4} \bigwedge^{3} V_5 \to 0 $
as well.
We will thus consider $W$ as a subspace of $\bigwedge^2 V_5$
and we identify $\bigwedge^5 V_5$ with $\IC$ through the volume form on $V_6$.
Then $W=\lambda_3(A)$.
As for the quadric hypersurface,
we call $Q(x)$ the projective quadric hypersurface associated to $q(x)\in \Sym^2W^\vee$. We have $Z=M_Z\cap Q(x)$. 
Consider $ w\in W$: since $W=\lambda_3(A)$, there is $\eta\in \bigwedge^3 V_5$ such that $\eta+x\wedge w\in A$.
By definition of $q: V_6\to \Sym^2 W^\vee$ given in \Cref{GM e EPW data}, the element $w=\lambda_3(\eta+x\wedge w)$ lies in $Q(x)$ if and only if $\eta\wedge w=\lambda_4(x\wedge \eta)\wedge w=0$.
\end{rmk}

\begin{rmk}
\label{secondo remark}
\normalfont
Consider instead some ordinary $GM$ intersection $Z=\IP(W)\cap G(2,V_5)\cap Q$, with associated $GM$ data $(W,V_6,V_5,\mu,q,\epsilon)$.
We consider again $W\subseteq \bigwedge^2 V_5\otimes L^\vee$, and $x\in V_6- V_5$ such that $Q=Q(x)$. As above $\bigwedge^3 V_6 \cong \bigwedge^3 V_5\oplus F_x$ and \eqref{sequenza spazi vettoriali} holds. Also $\bigwedge^3 V_5\cong (\bigwedge^2 V_5)^\vee$ via $\epsilon:\bigwedge^5 V_5\to L^{\otimes 2}$
and $A(Z)$ is
\[ \left\lbrace 
 (\eta,w)\in  \bigwedge^3 V_5\oplus W  \mbox{ } |\mbox{ } w\in W ,\mbox{ } 
\epsilon(\eta\wedge \mu(\_  ))_{|_W}= -(q(x)(w,-)\otimes [x]) \mbox{ as elements of }W^\vee 
\right\rbrace .\]

Now we forget $\epsilon$, we consider again $W\subseteq \bigwedge^2 V_5$ and $\bigwedge^3 V_5$ as $(\bigwedge^2 V_5)^\vee$.
We have then a splitting sequence 
$0\to  W^\perp \to  (\bigwedge^2 V_5)^\vee   \to  W^\vee    \to 0$.
Putting all together we can rewrite $A(Z)\subset 
W^\perp \oplus W^\vee\oplus F_x$ as
\begin{equation}\label{lagrangiana in tre pezzi}
 \left\lbrace 
 (\xi,-q(x)(w,-),x\wedge w)\in  W^\perp \oplus W^\vee\oplus  F_x \quad |\quad \xi\in W^\perp,\mbox{ }
 w\in W  
\right\rbrace .
\end{equation}
From this description, under the identification $\bigwedge^3 V_5\cong (\bigwedge^2 V_5)^\vee$ we have $A(Z)\cap \bigwedge^3 V_5=W^\perp$ by \cite[Proposition $3.13$, Item a)]{debarre2015}.\\
Consider the linear map $x\wedge W\to W^\vee$ that sends $x\wedge w$ to $-q(x)(w,-)$: its graph $\Gamma_x\subset W^\vee\oplus W$ induces a decomposition
$A(Z)=W^\perp\oplus \Gamma_x\subset W^\perp \oplus W^\vee\oplus F_x$.
\end{rmk}

\section{$K3$ surfaces whose associated double $EPW$ sextic is smooth}
\label{mio divisoriale}

We denote by $\mathcal{K}_{10}$ the moduli space of $\langle 10 \rangle$-polarized $K3$ surfaces.

\begin{defi} 
\label{definizione divisori per modulo di K3}
The divisor $D_{x,y}\subset \mathcal{K}_{10}$ is the locus of pairs $(S,H)$ such that there exists a primitive sublattice $\mathbb{Z}H+\mathbb{Z}D\subseteq \ns(S)$ whose Gram matrix is  
$\left[\begin{array}{cc}
10 & x  \\
x & y  \\ 
\end{array}\right].$
\end{defi}

In this section we want to find precise conditions on a Brill-Noether general, $\langle 10 \rangle$-polarized, $K3$ surface such that the associated double $EPW$ sextic $X_{A(S)}$ is a hyper-K\"ahler manifold.

We answer this question in two ways, as a condition on curves on the embedded $K3$ surface (\Cref{mio risultato su dominio: cor}) and as a divisorial condition on $\mathcal{K}_{10}$ (\Cref{mio risultato su dominio: versione pic}). 
In particular, thanks to \Cref{mio risultato su dominio: versione pic}, finding whether $X_A(S)$ is smooth becomes a lattice-theoretic problem on $\ns(S)$.

\begin{teo}\label{mio risultato su dominio: cor}
Let $S=\IP(W)\cap G(2,V_5)\cap Q$ be a $\langle 10\rangle$-polarized $K3$ surface. The double cover $X_{A(S)}$ of the associated $EPW$ sextic $Y_{A(S)}$ is smooth if and only if $S$ is strongly smooth and contains neither lines nor quintic elliptic pencils.
\end{teo}

We recall that $(S,H)\in \mathcal{K}_{10}$ is Brill-Noether general if and only if
the projective model of $(S,H)$ is a smooth ordinary Gushel-Mukai variety of dimension two.

From now on, we see a Brill-Noether general $(S,H)\in \mathcal{K}_{10}$ as a smooth ordinary Gushel-Mukai $2$-dimensional variety $S = \mathbb{P}(W) \cap G(2,V_5)  \cap Q$,
with Grassmannian hull $M_S=\mathbb{P}(W) \cap G(2,V_5)$. 

\begin{teo}
\label{mio risultato su dominio: versione pic}
Let $(S,H)\in \mathcal{K}_{10}$ be a $\langle 10 \rangle$-polarized $K3$ surface.
If $(S,H)$ is Brill-Noether general, let $A(S)\in \LG$ be the Lagrangian subspace associated to $(S,H)$. Then
\begin{itemize}
\item[1)] $(S,H)$ is Brill-Noether general if and only if $(S,H)\notin D_{h,0}$ for $h\in \lbrace 1,2,3 \rbrace$;
\item[2)] if $(S,H)$ is Brill-Noether general, then $(S,H)$ is strongly smooth if and only if 
$A(S)\notin \Sigma$, if and only if $(S,H)\notin D_{4,0}$.
\end{itemize}
Moreover, if $(S,H)$ is strongly smooth, then $A(S)\in \LG^0$ 
if and only if $(S,H)\notin D_{x,y}$, with $(x,y)\in \lbrace(1,-2),(5,0)\rbrace$.
In particular
\begin{itemize}
\item[3)] $Y^3_{A(S)}\cap \IP(V_5)=\emptyset$ if and only if $(S,H)\notin D_{1,-2}$;
\item[4)] $Y^3_{A(S)} - \IP(V_5)=\emptyset$ if and only if $(S,H)\notin D_{5,0}$.
\end{itemize}
\end{teo}

Grassmannian hulls $M_S$ of strongly smooth Brill-Noether general $K3$ surfaces are isomorphic, see \cite{iskovskikh1999} or \cite[Proposition $5.2$, Item $3)$]{ogrady2010}.
We need a characterization of the lines inside $M_S$. 
For $v\in V_5-\lbrace 0 \rbrace$ and $v\in V_3\subset V_5$ with $V_3\cong \mathbb{C}^3$, we set
$L_{v,V_3}=\lbrace  v\wedge t \mbox{ }|\mbox{ }t\in V_3 \rbrace$.
Every line in $G(2,V_5)$, hence in $M_S$, is of the form $\IP(L_{v,V_3})$.

\begin{lemma}
\label{rette in grassmannian hull}
Consider $v\in V_5-\lbrace 0 \rbrace$. There exists a three dimensional vector space $V_3\ni v$
such that the line $\IP(L_{v,V_3})\subset G(2,V_5)$ is contained in $M_S$
if and only if 
$\ker(P_v)\cap W$ has dimension exactly $2$.
In this case $L_{v,V_3}=\ker(P_v)\cap W$.
\begin{proof}
We fix a decomposition $V_5= \mathbb{C}v\oplus U$; we know that $\ker(P_v)=v\wedge U$, 
hence $L_{v,V_3}\subset\ker(P_v)$ for any $V_3\ni v$.
If $\ker(P_v)\cap W$ has dimension $2$, its projectivization is a line on $M_S$ of the form $\IP(L_{v,V_3})$.
Conversely, if $\IP(L_{v,V_3})\subset M_S$ then $L_{v,V_3}\subset W\cap \ker(P_v)$, 
in particular $\dim(W\cap \ker(P_v))\geq 2$.
The dimension is exactly $2$, or $M_S$ would contain a plane, which is absurd since $M_S$ has Picard rank $1$.
\end{proof}
\end{lemma}

When a $GM$ variety $Z$ is strongly smooth, for $v\in V_5-\lbrace 0\rbrace$ the kernel of the Pfaffian quadric $(P_v)_{|_W}$ is $\ker(P_v)\cap W$ and its corank is at least $\dim(W)-6$. For $x\notin V_6$ the kernel of $q(x)$ is computed in \cite[Proposition $3.13$]{debarre2015}, Item $b)$.

\begin{proof}[Proof of \Cref{mio risultato su dominio: cor}]

The proof consists of three parts. We already know, from \Cref{corrispondenze}, that strongly smoothness is a necessary condition to have $X_{A(S)}$ smooth. We need to prove $A(S)\notin \Delta$ i.e. $Y^3_{A(S)}=\emptyset$.

\medskip

\noindent \textit{$S$ contains no line if and only if $Y^3_{A(S)}\cap \mathbb{P}(V_5)=\emptyset$.}
Fix $x\in V_6-V_5$, so that $S=M_S\cap Q(x)$, with $Q(x)$ the projective quadric associated to $q(x)$; as in \Cref{primo remark} we have $\bigwedge^3 V_6\cong (\bigwedge^2 V_5)^\vee\oplus F_x$ and a linear map $\lambda_3 \colon (\bigwedge^2 V_5)^\vee\oplus F_x  \to  \bigwedge^2 V_5$ sending $(\phi, x\wedge w)$ to $w$.\\
For $v\in V_5-\lbrace 0 \rbrace$ we fix a decomposition $V_5=\IC v \oplus  U$.
Since $F_v=(v\wedge \bigwedge^2 U) \oplus v\wedge U$, we can write
\[ \lambda_3(A)=W \quad\quad \lambda_3(F_v)=v\wedge U \quad\quad  \ker((P_v)_{|_W})=(v\wedge U)\cap W  . \]
We have $A\cap (\bigwedge^2 V_5)^\vee=W^\perp$ by \Cref{secondo remark}; 
we can also write
$F_v\cap (\bigwedge^2 V_5)^\vee =v\wedge \bigwedge^2 U=\ker(P_v)^\perp$
and by \eqref{sequenza spazi vettoriali},
\begin{equation}
\label{conto dimensioni}
   \dim(A\cap F_v)=\dim (\lambda_3 (A\cap F_v))+\dim(A\cap F_v\cap (\bigwedge^2 V_5)^\vee).        
\end{equation}
We have $\lambda_3 (A\cap F_v)\subseteq \ker(P_v) \cap W$.
The dimension of $\ker(P_v) \cap W$ is at most $2$, see \Cref{rette in grassmannian hull}, in particular:
\begin{itemize}
\item[1)] If $\dim(\ker(P_v) \cap W)=2$, then $\ker(P_v)+W$ has dimension $9$, hence there is only one hyperplane containing both of them, so $\dim(W^\perp\cap \ker(P_v)^\perp)=1$;
\item[2)] If $\dim(\ker(P_v) \cap W)=1$, then $\ker(P_v)+W=\bigwedge^2 V_5$
and $W^\perp\cap \ker(P_v)^\perp=0$.
\end{itemize}
By \eqref{conto dimensioni}, since
$A\cap F_v\cap (\bigwedge^2 V_5)^\vee= W^\perp\cap \ker(P_v)^\perp$ in order to obtain $\dim(A\cap F_v)=3$ we must consider all the $v\in V_5$ for which 1) above holds.
For such a $v$ we have $\lambda_3(A\cap F_v)=\lambda_3(A)\cap \lambda_3(F_v)$, so that the RHS of \eqref{conto dimensioni} can be $3$.
The line $\ell=\IP(\lambda_3(A\cap F_v))$ lies in $M_S$, so we only need to prove $\ell\subset Q(x)$. 
By \Cref{primo remark}, we know that $q(x)(w,w)=0$ if and only if there there exists $\eta\in \bigwedge^3 V_5$ such that $\eta+x\wedge w\in A$ and $\eta\wedge w=0$.
Consider any $v\wedge u$ inside $\lambda_3(A\cap F_v)$: by hypothesis, there exists $\eta\in \bigwedge^3 V_5$ such that 
$\eta+x\wedge v\wedge u\in A\cap F_v$. This means that $\eta\in F_v$, so $\eta\wedge v\wedge u=0$. In particular $v\wedge u\in Q(x)$, hence the whole line $\IP(\lambda_3(A\cap F_v))$ is contained in $S$.

\smallskip

On the other hand, if $\IP(L_{v,V_3})\subset S$ for some $v\in V_5$, by \Cref{rette in grassmannian hull} we also have $L_{v,V_3}=\ker(P_v)\cap W\cong \IC^2$,
hence there is some non-zero $\eta\in W^\perp\cap \ker(P_v)^\perp$. Then to prove $[v]\in Y^3_A\cap \IP(V_5)$ we need $\lambda_3(A\cap F_v)=\lambda_3(A)\cap \lambda_3(F_v)=W\cap\ker(P_v)$.

Consider $v\wedge u\in W\cap\ker(P_v)$. We use \eqref{lagrangiana in tre pezzi}: 
under the decomposition $(\bigwedge^2 V_5)^\vee \cong W^\perp\oplus W^\vee$, the element
\[\alpha=(0,-q(x)(v\wedge u,-),x\wedge v\wedge u)\in W^\perp\oplus W^\vee \oplus F_x\]
lies in $A$ and $\lambda_3(\alpha)=v\wedge u$. To prove that $\alpha\in F_v$, so that we have $v\wedge u\in \lambda_3(A\cap F_v)$, we show $\alpha\wedge \eta=0$ for every $\eta\in F_v$, then we can conclude, since $F_v$ is Lagrangian. 
As $\alpha$ is zero on the first component, 
we are left to prove $0=\alpha\wedge \eta=-q(x)(v\wedge u,\eta)$ for $\eta\in F_v\cap W=L_{v,V_3}$. 
By hypothesis $\IP(L_{v,V_3})\subset S\subset Q(x)$:
this happens if and only if $q(x)_{|_{L_{v,V_3}\times L_{v,V_3}}}=0$ and in particular $q(x)(v\wedge u,-)=0$ on $L_{v,V_3}$.

\medskip

\noindent \textit{If $Y^3_{A(S)}-\IP(V_5)\neq \emptyset$, then $S$ contains a quintic elliptic pencil.}
Consider $[x]\in Y^3_A - \IP(V_5)$, so that $S$ is the transverse intersection of $M_S$ and $Q(x)$: by \cite[Proposition $3.13$, Item b)]{debarre2015} the corank of $q(x)$ is $3$.\\
The quadric $Q(x)$ is a cone of $\IP(\ker (q(x)))=\IP(A\cap F_x)$ over a smooth quadric surface $\mathcal{Q}\subset\IP^3$; the latter has two families of lines on it and two lines inside $\mathcal{Q}$ intersect if and only if they do not lie in the same family. Let $\lbrace\ell_t\rbrace_{t\in \IP^1}$ be one of the two families. The cone of $\IP(\ker (q(x)))$ over $\ell_t$ is a $\IP^4$ which we denote by $\pi_t$.
By transversality, the intersection $M_S \cap \pi_t = S\cap \pi_t$ is a degree-$5$ curve on $S$.\\
Consider now $(\pi_t\cap S)$ and $(\pi_{t'}\cap S)$ for $t\neq t'$. We have
$\pi_t\cap \pi_{t'}=\IP(\ker (q(x)))$, as $\ell_t$ and $\ell_{t'}$ are in the same family, so $(\pi_t\cap S)\cap (\pi_{t'}\cap S)=\IP(\ker (q(x)))\cap S$. But the latter is empty, since $\IP(\ker (q(x)))=\sing(Q(x))$ and $S$ is smooth. So the pencil $\lbrace S\cap \pi_t\rbrace_{t\in \IP^1}$ is elliptic.

\medskip

\noindent \textit{If $S$ contains a quintic elliptic pencil, then $Y^3_{A(S)}-\IP(V_5)\neq \emptyset$.}
We call $E\in \ns(S)$ the class of the elliptic pencil: 
the sublattice $\IZ H+\IZ E=K\subset \NS(S)$ has Gram matrix 
$\begin{pmatrix}
10 & 5\\
5 & 0
\end{pmatrix}$. The orthogonal complement of $H$ in $K$ is generated by $\kappa=H-2E$; we denote by $\kappa_2$ the corresponding class in $\NS(S^{[2]})$.

From now on, we consider notations and results from \cref{subsec involuzione}.
Following that, we call $e,f$ two canonical generators of a copy of $U$, moreover $e_1,f_1$ will be the canonical generators of a second copy of $U$, orthogonal to the first one.
By Eichler Criterion \cite[Proposition 3.3]{gritsenko2008}, there is an isometry $\psi:H^2(S^{[2]},\IZ)\to \Lambda$ such that 
$\psi(H_2-2\delta)=h=e+f$ and $\psi(\kappa_2)=g+2\ell-2e_1-8f_1$;
the involution $j$ on the moduli space exchanges $g$ and $\ell$ in $\Lambda_h$, so that $j(\kappa_2)=g+2\ell-2e_1-8f_1$ is an algebraic class for $X_{A(S)}$, which lies in $D_A^\perp$.
It has square $10$ and divisibility $2$ in $H^2(X_{A(S)},\IZ)$, so $D_A$ lies in a flopping wall and is not ample, see \cite[Theorem $5.1$, Item $b)$]{debarre2017}.

Hence $Y_{A(S)}^3\neq \emptyset$ and more precisely $Y_{A(S)}^3 - \IP(V_5)\neq \emptyset$, since points in $Y_{A(S)}^3 \cap \IP(V_5)$ correspond to lines in $S$, hence sublattices in $\ns(S)$ whose Gram matrix is
$\begin{pmatrix}
10 & 1\\
1 & -2
\end{pmatrix}$.
It is then sufficent to observe that such a sublattice is neither an overlattice nor a sublattice of $H\IZ+E\IZ$.
\end{proof}

We can now prove the second version of the characterization.

\begin{proof}[Proof of \Cref{mio risultato su dominio: versione pic}]
For $1)$ see \cite[Lemma $2.8$]{greer2014}. We do not need to ask $(S,H)\notin D_{5,2}$, as in \cite{greer2014}, since this is a divisor inside the complement of $\mathcal{K}_{10}$ in its closure; to see that, one checks that the orthogonal complement of $H$ in $\mathbb{Z} H+\mathbb{Z} D$ has square $-2$ and concludes by \cite[Theorem 2.7]{gritsenko2010}.
For $2)$ see \cite[Lemma $2.7$]{greer2014}.
The rest of the statement comes directly from \Cref{mio risultato su dominio: cor}.
Lines on $S=\IP(W)\cap G2,V_5)\cap Q$ correspond to sublattices of $\ns(S)$ whose Gram matrix is
$\left[\begin{array}{cc}
10 & 1 \\
1 & -2 
\end{array}\right]$,
To prove $3)$ it is then sufficient to show that this lattice admits no non-trivial overlattice.
For $4)$ the proof is the same, with the lattice whose Gram matrix is
$\left[\begin{array}{cc}
10 & 5 \\
5 & 0 
\end{array}\right]$.
\end{proof}

Finally, the description of lines in $S$ can be made more precise.

\begin{prop}
\label{rette su S}
Let $S= \mathbb{P}(W) \cap G(2,V_5)  \cap Q$ be a strongly smooth $K3$ surface.
Let $A(S)$ be the associated Lagrangian subspace. Then
\[ |Y^3_{A(S)}\cap \mathbb{P}(V_5)|= |\lbrace \ell \subset S| \ell \mbox{ line}\rbrace|. \]
\end{prop}
\begin{proof}
This is an easy consequence of the first part of the proof above. To show the one-to-one correspondence,
one observes that two lines $\IP(L_{v,V_3})$ and $\IP(L_{v',V'_3})$ intersect in at most one point if $[v]\neq [v']$.
\end{proof}

\section{Inducing automorphisms on double EPW sextics}
\label{auto indotti}

We explain here how we can use our main result, and particularly Theorem \ref{mio risultato su dominio: versione pic}, to deduce the existence of automorphisms on families of smooth double $EPW$ sextics. This was actually the first motivation for our study of strongly smooth $GM$ surfaces, since the automorphisms of $S$ act on the associated $EPW$ sextic.

\begin{prop}(Debarre and Kuznetsov)\label{auto indotti su A(Z)}
For $Z$ a strongly smooth $GM$ variety, there is a natural inclusion of $\Aut(Z,\IP(W))$ in $\Aut(Y_{A(Z)})\cong \Aut(Y_{A(Z)^\perp})$ as the stabilizer of $[V_5]\in Y_{A(Z)^\perp}$.
\end{prop}
\begin{proof}
See \cite[Proposition $3.21$, Item $c)$]{debarre2015}: the result still holds if we replace "dimension $\geq 3$" in the statement with "$X$ strongly smooth".
\end{proof} 

O'Grady showed that, for $S=\IP(W)\cap G(2,V_5)\cap Q$ very general, the associated double $EPW$ is smooth, however the group $\Aut(S,\IP(W))$ is trivial for $S$ very general! In addition, all those automorphisms turn out to be symplectic, see \cref{tuuuto symply} and \cref{gli auto di una GM che e una K3}, and families of $K3$ surfaces in $\mathcal{K}_{10}$ carrying symplectic actions typically have big codimension.

So, to find automorphisms on some $S$ that induce automorphisms on a \emph{smooth} double $EPW$ sextic, the tricky part is to control the generality. 
By \Cref{mio risultato su dominio: versione pic}, we can easily deduce the existence of automorphisms on some smooth double $EPW$ sextic, appealing to the classification of automorphisms on $K3$ surfaces.
Let $(V_6,V_5,A)$ be the Lagrangian data associated to $S$: by \Cref{GM e EPW data}, $[V_5]\in Y^3_{A^\perp}$, hence $A$ always lies in $\Pi$.

In the following statement, we denote by $D_n$ the dihedral group of order $2n$.
In the proof we use some results from the next section; we decided to postpone them to keep the exposition more cohesive.

\begin{prop}
\label{famiglie defo di K3 GM simply}
\label{auto indotti da K3 su EPW}
Let $G$ be one of the following groups,
\[
\IZ/n\IZ \mbox{ for }n\in\lbrace 2,3,4 \rbrace,\quad (\IZ/2\IZ)^2, \quad D_n\mbox{ for }n\in\lbrace 4,5,6\rbrace.
\]
There is a family of Lagrangian subspaces in $\Pi - (\Sigma \cup \Delta)$ such that, for any $A$ in the family, the associated double $EPW$ sextic $X_{A}$ 
is smooth and admits a symplectic action of $G$ which commutes with the covering involution.\\
For a general element of the family, $\Aut_{D_{A}}(X_{A})=\iota_{A}\times G$.
\begin{proof}

A $\langle 10 \rangle$-polarized $K3$ surface $(S,H)$ admits a symplectic action of $G$ which fixes $H$ if and only if its Néron-Severi group contains $\mathcal{L}_G=\IZ H\oplus \Omega_G$, where $\Omega_G$ is a negative-definite lattice associated to $G$, which turns out to be the coinvariant lattice for the action of $G$ in cohomology.
This has been proved by Nikulin \cite{nikulin1979a}, but see also \cite[Proposition 6.3]{garbagnati2009}.

We start by considering $G=\IZ/n\IZ$ with $n=2,\ldots 6$, $(\IZ/2\IZ)^2$ or $D_4$ and we want to prove that there is a family of $K3$ surfaces for which the very general element $S$ has $\NS(S)=\mathcal{L}_G$. 
Then we will have $G=\Aut(S,\IP(W))$ for $G\neq \IZ/n\IZ$ with $n=5,6$, and $\Aut(S,\IP(W))=D_n$ otherwise: in the latter case we appeal to \cite[Propositions 8.1 and 9.1]{garbagnati2013} and then \Cref{gli auto di una GM che e una K3}, in the former simply to \Cref{gli auto di una GM che e una K3}.

By the surjectivity of the period map for marked $K3$ surfaces, proving the existence of the family amounts to showing that there is a primitive embedding $\mathcal{L}_G$ in $\Lambda_{K3}=U^{\oplus 3}\oplus E_8(-1)^{\oplus 2}$.
By \cite[Propositions 6.2 and 7.7]{garbagnati2013}, there is a family of elliptic $K3$ surfaces whose Néron-Severi group is $U\oplus \Omega_G$. Now, any embedding of $U$ in an even lattice $M$ induces a splitting $M=U\oplus U^\perp$, in particular in our case we have $\Omega_G\subset U^{\oplus 2}\oplus E_8(-1)^{\oplus 2}\subset \Lambda_{K3}$. We fix a canonical basis $\lbrace e,f \rbrace$ for the distinguished copy of $U$ and we consider the lattice $\IZ(e+5f)\oplus \Omega_G$, which is a copy of $\mathcal{L}_G$ and is primitive in $\Lambda_{K3}$, since the two summands lies in two different direct summands of the big lattice.
This give us a maximal family of $K3$ surfaces with a symplectic action of $G$.

For the very general element $(S,H)$ in the family, the $\langle 10 \rangle$-polarization has divisibility 10 in the Néron-Severi group, in particular $(S,H)$ is Brill-Noether general and strongly smooth by \Cref{mio risultato su dominio: versione pic}. Since these properties fail on a finite number of divisors in the moduli space, this means that also the \textit{general} element in the family enjoys them (these families has positive dimension).
For the same lattice-theoretic reasons, $(S,H)$ will not lie in $D_{h,0}$ for $h=4,5$ nor $D_{1,-2}$, so its associated double $EPW$ sextic is smooth and carries a symplectic action of $G$ fixing the polarization $D_A$, by \Cref{auto indotti su A(Z)} and \Cref{tuuuto symply}. 

\Cref{mia sequenze con auto YA versione K3} ensures that the group of automorphisms fixing $D_A$ on the associated $EPW$ sectic is isomorphic to $\Aut(S,\IP(W))$ for a general $(S,H)$, since in that case the projective model does not contain any line or conic. 
\end{proof}
\end{prop}

\section{Bounds on automorphisms of Gushel-Mukai varieties}
\label{automorfismi vari e eventuali}

\subsection{Bounds for Brill-Noether general K3 surfaces}
\label{auto di K3 che sono GM}

In the next two sections we provide some bounds on the automorphism group of $GM$ varieties.
Grassmannian hulls will play a fundamental role and the key point is the following observation.

\begin{lemma}
\label{auto di gm vengono da hull per K3}
Let $Z$ be a strongly smooth $GM$ variety with Grassmannian hull $M_Z$. Then $\Aut(Z,\IP(W))\subset \Aut(M_Z)$.
\end{lemma}
\begin{proof}
Since $S$ spans $\IP(W)$, the restriction $\rho\colon\lbrace \alpha \in \Aut(\IP(W)) \mbox{ } | \mbox{ }\alpha(S)=S\rbrace\to \Aut(S,\IP(W))$ is an isomorphism of groups.
By \cite[Corollary $2.11$]{debarre2015}, any $\alpha\in \Aut(S,\IP(W))$ is induced by an automorphism $(g_W,g_{V},g_L)$ of ordinary $GM$ data, where the isomorphism $[g_W]\in \PGL(W\otimes L)=\PGL(W)$ is the restriction of $\left[\bigwedge^2 \left({g_V}_{|_{V_5}}\right)\right]$ to $\IP(W)$.
Hence $\rho$ factors as
$\lbrace \alpha \in \Aut(\IP(W)) \mbox{ } | \mbox{ }\alpha(S)=S\rbrace\xrightarrow{\rho_1} \Aut(M_S) \xrightarrow{\rho_2} \Aut(S,\IP(W))$.
Since $\rho$ is injective, the same holds for $\rho_1$.
\end{proof}

Let $S=\IP(W)\cap G(2,V_5)\cap Q$ be a Brill-Noether general $K3$ surface.
The automorphism group $\Aut(S,\IP(W))$ 
is the group of the automorphisms of the abstract surface $S$ fixing $H$; this group is finite, see \cite[Chapter $5$, Proposition $3.3$]{huybrechts2016}.

O'Grady \cite{ogrady2010} studied in detail the 6-dimensional complete linear system $|H_2-2\delta|$. 
This linear system is naturally isomorphic to the space of quadrics in $\IP(W)$ containing $S$ and the associated map admits a very explicit description (\cite[$(4.2.5)$]{ogrady2010}),
\begin{align}\label{definizione di g}
\phi\colon S^{[2]} & \dashrightarrow |H_2-2\delta|^\vee\cong \IP^5 \\
Z & \mapsto \lbrace  \mathcal{Q}\in |V_6| \mbox{ such that } \langle Z \rangle\subset \mathcal{Q}  \rbrace .
\end{align}
The map $\phi$ factors as a small contraction
$c\colon S^{[2]}\to X_{A(S)^\perp}$ and $f\colon X_{A(S)^\perp}\dashrightarrow \IP^5$. The latter can be identified to the double cover $f_{A(S)\perp}$ up to a finite number of flops \cite[Theorem 4.15]{ogrady2010}, in particular it is generically 2-to-1 and induces a birational involution $\iota$ on $S^{[2]}$ which acts as minus the reflection with respect to the class $H_2-2\delta$ in cohomology \cite[Proposition 4.20]{ogrady2010}.
The birational transformation $c$ induces an identification of $\Aut_{D_{A(S)^\perp}}(X_{A(S)^\perp})$ with $\Bir_{H_2-2\delta}(S^{[2]})$, the group of birational endomorphisms fixing $H_2-2\delta$ and \eqref{sequenza lift alla doubleEPW} translates as
\[  1\to \lbrace \id,\iota \rbrace\to \Bir_{H_2-2\delta}(S^{[2]}) \xrightarrow{\lambda} \Aut(Y_{A(S)^\perp}) \to 1. \]

We denote by $\mathfrak{G}_n$ the symmetric group on $n$ objects and by $\mathfrak{A}_n$ the alternating group on $n$ objects.
\begin{prop}
\label{gli auto di una GM che e una K3}
If $S$ is strongly smooth, $\Aut(S,\IP(W))$ acts symplectically on $S$ and is isomorphic to one of the following groups:
\begin{equation}\label{possibili gruppi per dim 2}
\begin{gathered}
   \IZ/n\IZ \mbox{ for }n= 1,2,3,4, ,\quad \quad  D_n \mbox{ for } n= 2,\ldots,6    \\
          \mathfrak{A}_4, \quad \quad \quad   \mathfrak{G}_4,\quad \quad \quad  \mathfrak{A}_5.\quad \quad 
\end{gathered}
\end{equation}
\begin{proof}

We prove simplecticity.
We know that \eqref{sequenza lift alla doubleEPW} splits, identifying the lifting of $\Aut(Y_{A(S)^\perp})$, which we denote by $\mathcal{S}$, with the group of automorphisms acting trivially on $H^{0,2}(X_{A(S)^\perp})$.
Since $c$ is an isomorphism in codimension one, $\mathcal{S}$ seen as a subgroup of $\Bir_{H_2-2\delta}(S^{[2]})$ acts simplectically on $S^{[2]}$. Inside $\mathcal{S}$ there is $G$, the lifting of $\Aut(S,\IP(W))\subset \Aut(Y_{A(S)^\perp})$. We denote by $G'\subset \Bir_{H_2-2\delta}(S^{[2]})$ the group of natural automorphisms $\alpha^{[2]}$ on $S^{[2]}$ induced by $\alpha\in \Aut(S,\IP(W))$: if we can prove $G=G'$ we are done, since $\alpha^{[2]}$ is symplectic if and only if $\alpha$ is.\\
Fix some $\alpha\in \Aut(S,\IP(W))\subset \Aut(Y_{A(S)^\perp})$. A straightforward computation, using \eqref{definizione di g} and \eqref{diagrammi auto GM data}, gives $\phi\circ\alpha^{[2]}=\alpha\circ \phi$, hence $\lambda(G)=\lambda(G')$.
For $\tilde{\alpha}\in G$ the lifting of $\alpha$, we have then $\tilde{\alpha}=\alpha^{[2]}$ or $\tilde{\alpha}=\alpha^{[2]}\circ \iota$. Note that $\alpha^{[2]}$ fixes the class $H_2$, while $\iota$ only fixes $H_2-2\delta$ and its multiples \cite[Proposition 4.21 (b)]{ogrady2005}, so the first equality holds if and only if $\tilde{\alpha}$ fixes $H_2$.\\
We turn to $|H_2|$, as a linear system on $X_{A(S)^\perp}$, via the small contraction $c$; it is the pullback via the Hilbert-Chow contraction of $H^{\otimes 2}$ on the symmetric product $S^{(2)}$ , 
\[    H^0(S^{[2]},H_2)=H^0(S^{(2)},H^{\otimes 2})=\Sym^2 W^\vee\]
The action of the lifting $\tilde{\alpha}$ of $\alpha$ to $X_{A(S)^\perp}$ is induced by the action of $SL(V_6)$ \cite[Proof of Proposition A.2]{debarre2020}, thus $\tilde{\alpha}$ acts on the global sections of $H_2$ since it induces a linear automorphism of $W$, see \Cref{diagrammi auto GM data}. So $\tilde{\alpha}^*H_2=H_2$, which provides us $\tilde{\alpha}=\alpha^{[2]}$ as desired.

To obtain the bound we observe that, by \Cref{auto di gm vengono da hull per K3}, $G$ is a finite subgroup of $\Aut(M_S)$ and the latter is isomorphic to $\PGL(2,\IC)$ (\cite[Theorem $7.5$]{piontkowski1999}), since the strongly smoothness of $S$ implies the condition of generality for $M_S$ in \cite{piontkowski1999}, cf. \cite[Proposition $2.22$]{debarre2015}.
So $G$ is either cyclic, dihedral, $\mathfrak{A}_4$, $\mathfrak{G}_4$ or $\mathfrak{A}_5$. For cyclic groups the order is at most $8$ by \cite[Theorem $4.5$]{nikulin1979a} and actually when the order is $n=7,8$ the lattice $\Omega_G^\perp$ (same notation as in \Cref{auto indotti da K3 su EPW}) does not contain an element of square $10$, as a computation modulo $2n$ shows, so we can rule them out.
Moreover, if $G$ contains an order $n=5,6$ element, then $D_n\subset G$ by \cite[Proposition 8.1]{garbagnati2013}. Finally, $n\leq 6$ for $D_n$ by \cite{xiao1996a}. 
\end{proof}
\end{prop}

\begin{rmk}
\normalfont
In \Cref{famiglie defo di K3 GM simply}, we provide some families of $K3$ surfaces 
with prescribed $\Aut(S,\IP(W))$. According to \Cref{gli auto di una GM che e una K3}, the possible cases left are $\Aut(S,\IP(W))\in \lbrace D_3,\mathfrak{A}_4,\mathfrak{G}_4,\mathfrak{A}_5 \rbrace$.
\end{rmk}

When $S$ contains no line, $\phi:S^{[2]}\dashrightarrow |H_2-2\delta|$ behaves particularly well.
Let $\mathcal{C}$ be the set of smooth conics inside $S$: by \cite[Claim 4.19]{ogrady2010} $|H_2-2\delta|$ is base-point-free and 
\begin{equation}\label{sing e coniche}
Y^3_{A(S)^\perp}=\bigcup_{C\in \mathcal{C}} \phi(C^{(2)}).
\end{equation}

When moreover $S$ is strongly smooth, we are able to provide a bound for the automorphism group of the associated $EPW$ sextic $Y_{A(S)}$, which is particularly useful when there is no conic on $S$; we used \Cref{mia sequenze con auto YA versione K3} under this condition in the proof of \Cref{famiglie defo di K3 GM simply}, in the previous section.

\begin{defi}
We denote by $\Aut_\mathcal{C}(S,\IP(W))\leq \Aut(S,\IP(W))$ the subgroup of automorphisms acting trivially on $\mathcal{C}$.
\end{defi}
\noindent Clearly $\Aut_\mathcal{C}(S,\IP(W))$ appears in \eqref{possibili gruppi per dim 2}.

\begin{prop}
\label{mia sequenze con auto YA versione K3}
Let $S=\IP(W)\cap G(2,V_5)\cap Q$ be a strongly smooth Brill-Noether general $K3$ surface with $N$ smooth conics on it. 
Suppose that $S$ contains no line: the automorphism group of $Y_{A(S)}$ sits in an exact sequence 
\begin{equation}\label{SEC auto YA quando ce K3}
1\to \Aut_\mathcal{C}(S,\IP(W)) \to \Aut(Y_{A(S)}) \xrightarrow{g} \mathfrak{G}_{N+1}
\end{equation}
where $g$ sends $\alpha\in \Aut(Y_{A(S)})$ to its action on $Y^3_{{A(S)}^\perp}$, whose cardinality is $N+1$.
\begin{proof}
The morphism $g$ in \eqref{SEC auto YA quando ce K3} is well defined by \Cref{auto di EPW e suo duale sono uguali} and its kernel is the subgroup of automorphisms of $Y_A$ fixing pointwise $Y^3_{A^\perp}$. 
We are only left to prove that $\Aut_\mathcal{C}(S,\IP(W))=\ker(g)$, we actually prove a bit more. 
Indeed, by \eqref{sing e coniche} the fibers over $Y^3_{A(S)^\perp}$ are given by $P_S$ and $C^{(2)}$ for $C\in \mathcal{C}$; by the same computations of \Cref{gli auto di una GM che e una K3}, the action of $\alpha$ on $\mathcal{C}$ induces a permutation on $\lbrace P_S,C^{(2)}_1, \ldots , C^{(2)}_N \rbrace$ hence on $Y^3_{A(S)^\perp}$, which is the same as the one induced by $\alpha$ as an element of $\Aut(Y_{A(S)})$ (note that, since $\alpha$ acts on $M_S$ too by \Cref{auto di gm vengono da hull per K3}, $\alpha^{[2]}$ always fixes $P_S$).  
\end{proof}
\end{prop}

\begin{rmk}\label{card fiber}
\normalfont
The morphism $g$ in \eqref{SEC auto YA quando ce K3} is not necessarily surjective.
Indeed, let $(S,H)\in\mathcal{K}_{10}$ be a very general element of $D_{2,-2}$ (see \Cref{definizione divisori per modulo di K3}). Then $S$ contains exactly one smooth conic, say $C$.
An involution of $Y_{A(S)}$ exchanging the two elements inside $Y^3_{A(S)^\perp}$ would lift to a symplectic involution on $X_{A(S)}$ (which is smooth for $(S,H)$ general, see \Cref{mio risultato su dominio: versione pic}), but then $X_{A(S)}$ would have Picard rank at least $9$ by \cite[Corollary $5.2$]{mongardi2016arx}. This is impossible since $X_{A(S)}$ and $S^{[2]}$ have the same Picard rank, cf. \Cref{pic di EPW e Hilb}.

Furthermore, there exists $(S',H')\in D_{2,-2}$ such that $(S^{[2]}, H_2-2\delta)$ and $((S')^{[2]}, H'_2-2\delta')$ are birational as varieties with a big and nef line bundle, or equivalently $A(S)=A(S')$ and $P_{S'}$ is sent to $\phi(C^{(2)})$ via the map associated to $|H'-2\delta'|$. 
The argument above proves that $(S,H)$ and $(S',H')$ are not isomorphic.
\end{rmk}

\subsection{Bounds in greater dimension}

In \Cref{auto di K3 che sono GM} we dealt with $2$-dimensional $GM$ varieties. Now we consider $GM$ varieties of dimension $3$ and $4$.
By \cite[Proposition $3.21$]{debarre2015} their automorphism group is finite, and trivial in the general case.

\begin{rmk}
\normalfont
By \Cref{auto indotti su A(Z)}, whenever the associated Lagrangian $A(Z)$ does not lie in $\Delta\cap \Pi$, automorphisms of prime order of $Z$ have order at most $11$, because the automorphism lifts to $X_{A(Z)}$ and $X_{A(Z)^\perp}$, of which at least one is smooth, and a symplectic automorphism of prime order on a deformation of a Hilbert scheme of two points has order at most $11$, see \cite[Corollary 2.13]{mongardi2013a}. 
\end{rmk}

\begin{prop}
\label{gruppi se fissi punto in Y2A}
Let $Z$ be a smooth ordinary $GM$ variety of dimension $3$ and $G\subset \Aut(Z)$. Then there exists a group $H$ among
\begin{equation}\label{sequenza sopra caso Y2A}
1, \quad \quad  A_4, \quad \quad  \mathfrak{G}_4, \quad \quad  A_5, \quad \quad \IZ/k\IZ \quad\mbox{or} \quad   D_k \quad\mbox{for }k\geq 2 
\end{equation}
and $r\geq 1$ such that $G$ sits in an exact sequence $1 \to  \IZ/r\IZ   \to G \to H \to 1 $.
\end{prop}
\begin{proof}
The Grassmannian hull $M_Z$ is smooth since $Z$ is three-dimensional and $\Aut(Z,\IP(W))=\Aut(Z)$ since $Z$ has Picard rank one. By \cite[Theorem $6.6$]{piontkowski1999}, there is a short exact sequence
\begin{equation}
 1\to \mathbb{C}^4 \rtimes \mathbb{C}^* \to \Aut(M_Z)\xrightarrow{\pi} \PGL(2,\mathbb{C}) \to 1 ,
\end{equation}
thus $G$ is an extension of $H=\pi(G)$ by $N=G\cap (\mathbb{C}^4 \rtimes \mathbb{C}^*)$.
The group $G$ is finite by \Cref{auto di gm vengono da hull per K3}, so $H\subset \PGL(2,\mathbb{C})$ has to appear in \eqref{sequenza sopra caso Y2A}.
As for $N$, it is finite inside $\mathbb{C}^4 \rtimes \mathbb{C}^*$, so it is isomorphic to a finite subgroup of $\IC^*$.
\end{proof}

\begin{rmk}
\normalfont
In \cite{debarre2021}, Debarre and Mongardi produced an explicit example of $EPW$ sextic with automorphism group isomorphic to $\PSL(2,\mathbb{F}_{11})$. This allows them to find examples of smooth ordinary $GM$ threefolds with automorphism group isomorphic to $\IZ/11\IZ$, $D_6$, $\IZ/6\IZ$, $\IZ/3\IZ$, $D_5$, $D_2$, $\IZ/2\IZ$.
\end{rmk}

The \emph{spin group} $\Spin(n)$, for $n\geq 2$, is the universal cover of the special orthogonal group $\SO(n)$. 
It is possible to study finite subgroups of $\Spin(5)$ using this definition: finite subgroups of $\SO(5)$ have been classified, see \cite[Corollary $2$]{mecchia2011}.

The quotient by $\lbrace \pm 1\rbrace$ in the result below comes from the fact that, for any $\alpha\in \PGL(V_5)$ acting on the Grassmannian hull $M_Z$, a representative in $\GL(V_5)$ of the form
$\begin{pmatrix}
\bar{T} & \underline{0}\\
{^t\underline{b}} & b
\end{pmatrix}$
is chosen, such that $\dete(\bar{T})=1$. A representative in this form for $\alpha$ is unique, up to multiplication by $-1$. See \cite[Proposition $5.2$]{piontkowski1999} for details.

\begin{prop}
\label{gruppi se fissi punto in Y1A}
Let $Z$ be a smooth ordinary $GM$ variety of dimension $4$ and $G\subset \Aut(Z)$.
Then $G$ is the quotient by $\lbrace \pm 1\rbrace$ of a group $\tilde{G}$ such that there is an exact sequence
\[ 1\to   \IZ/r\IZ \to  \tilde{G} \to H \to 1     \]
for some $r\geq 1$ and $H$ a finite subgroup inside $\Spin(5)$.
\begin{proof}
As in \Cref{gruppi se fissi punto in Y2A}, we just have to impose finiteness for a subgroup of $\Aut(M_Z)$.
By \cite[Proposition $5.2$]{piontkowski1999}, the group $G$ is an extension of a subgroup of $\left( \Sp(4,\mathbb{C})\times \mathbb{C}^* \right)  /\lbrace \pm 1\rbrace$ by some $N\subset \mathbb{C}^4$.
Actually $N=\lbrace 1 \rbrace$ as $G$ is finite, hence $G$ is isomorphic to
the quotient of a finite subgroup of $\Sp(4,\mathbb{C})\times \mathbb{C}^*$.
Any finite subgroup of $\Sp(4,\mathbb{C})$ lies inside its maximal compact subgroup, which is the compact symplectic group $\Sp(2)$.
In turn $\Sp(2)$ is isomorphic to $\Spin(5)$.
\end{proof}
\end{prop}

\begin{rmk}
\normalfont
In \cite{debarre2021}, Debarre and Mongardi produced smooth ordinary $GM$ fourfolds with automorphism group isomorphic to $D_6$, $\IZ/3\IZ$, $D_5$, $\IZ/5\IZ$, $D_2$, $\IZ/2\IZ$.
\end{rmk}

\begin{cor}
Consider $A\in \LG - \Sigma$ and $G\subseteq \Aut(Y_A)$.
Suppose that $G$ fixes a point $[v]\in Y^k_A$ or a point $[V_5]\in Y^k_{A^\perp}$ for $k\in \lbrace1,2\rbrace$. 
\begin{itemize}
\item[1)] If $k=2$, then $G$ sits in an exact sequence $1 \to  \IZ/r\IZ   \to G \to H \to 1 $ for some $r\geq 1$ and $H$ a group in \eqref{gruppi se fissi punto in Y2A}. 
\item[2)] If $k=1$, then $G$ is the quotient by $\lbrace \pm 1\rbrace$ of a group $\tilde{G}$, where $\tilde{G}$ is an extension of some finite subgroup of $\Spin(5)$ by a cyclic group $\IZ/r\IZ$.
\end{itemize}
\end{cor}
\begin{proof}
We consider the case of $G$ fixing $[V_5]\in Y^k_{A^\perp}$, the dual case
follows by \Cref{auto di EPW e suo duale sono uguali}. 
We denote by $Z$ the ordinary $GM$ variety associated to the Lagrangian data $(V_6,V_5,A)$: by \Cref{strongly smooth e no decomponibili: cor} it is a strongly smooth $GM$ variety of dimension $5-k$.
By \Cref{auto indotti su A(Z)}, there is an inclusion $G\subset \Aut(Z,\IP(W))$, then the result follows from \Cref{gruppi se fissi punto in Y2A} and \Cref{gruppi se fissi punto in Y1A}.
\end{proof}

\noindent If the associated Lagrangian $A(Z)$ lies in $\Delta$ or $\Sigma$ we can be more precise.

\begin{defi}
A plane $\pi\subset G(2,V_5)$ is a \emph{$\sigma$-plane} if it is of the form $\IP(v\wedge V_4)$ for some $4$-dimensional vector space $V_4\subset V_5$ and $v\in V_4-\lbrace 0\rbrace$.
A \emph{$\tau$-quadric surface} in $G(2,V_5)$ is a linear section of $G(2,V_4)$ for some $4$-dimensional $V_4\subset V_5$. 
\end{defi}

\begin{cor}
\label{bound sugli auto se Y3A o Y3Ap nonvuoto}
Let $Z$ be a smooth ordinary $GM$ variety of dimension $4$ and associated Lagrangian subspace $A(Z)$. If one of the following holds:
\begin{itemize}
\item $Z$ contains a $\sigma$-plane,
\item $Z$ contains a $\tau$-quadric surface,
\end{itemize}
then there is an exact sequence $1\to N \to \Aut(Z) \to \mathfrak{G}_{n} $,
where $n=|Y^3_{A(Z)}|$ (resp. $n=|Y^3_{A(Z)^\perp}|$) if $Z$ contains a $\sigma$-plane (resp. a $\tau$-quadric surface),
and $N$ is one of the groups in the list at \eqref{possibili gruppi per dim 2}.
\begin{proof}
By \cite[Theorem $4.3$, Item $c)$]{debarre2019} for $\sigma$-planes, \cite[Remark $5.29$]{debarre2019} and \cite[Section $7.3$]{debarre2014} for $\tau$-quadric surfaces, the associated Lagrangian subspace $A(Z)$ lies in $\Delta$ or $\Pi$. So there is a strongly smooth Brill-Noether general $K3$ surface whose associated Lagrangian subspace is $A(Z)$. Then we use \Cref{mia sequenze con auto YA versione K3}.

\end{proof}
\end{cor}

The coarse moduli space of smooth $GM$ varieties of dimension $4$ is constructed in \cite{debarre2020a}. Inside it, the family of ordinary $GM$ varieties containing a $\sigma$-plane (resp. a $\tau$-quadric surface) has codimension $2$ (resp. $1$), see \cite[Remark $5.29$]{debarre2019};
the general members of both families are rational, see \cite[Proposition $7.1$ and Proposition $7.4$]{debarre2014}.

\appendix

\section{An involution of the moduli space}
\label{subsec involuzione}

We collect some results which are useful for the proof of \cref{mio risultato su dominio: cor}, but they are different in flavour from the rest of the paper;
for details about the construction in this section we refer to \cite[Section $3.2$]{debarre2017}.

We denote by $\Lambda=U^{\oplus 3}\oplus E_8(-1)^{\oplus 2}\oplus \langle -2 \rangle$ the abstract lattice isomorphic to $H^2(X,\IZ)$ for $X$ equivalent by deformation to the Hilbert square on a K3 surface. 
We fix a copy of the hyperbolic plane $U \subset \Lambda$, with canonical basis $\lbrace e,f\rbrace$.\\
We take $h=e+f$ and we call $\Lambda_h\subset \Lambda$ the orthogonal complement of $h$.
The group of isometries of $\Lambda$ fixing $h$ acts on $\Lambda_h$ as $\tilde{O}(\Lambda_h)=\lbrace \theta\in O(\Lambda_h) $ $|$ $\bar{\theta}=\id_{A_{\Lambda_h}} \rbrace$, cf. \cite[Proposition $3.12$, Item $i)$]{gritsenko2010a}.
We denote $\mathcal{Q}_2=\lbrace [x]\in\IP(\Lambda_h\otimes\IC)\mbox{ }|\mbox{ }x^2=0,\mbox{ }(x,\bar{x})>0\rbrace$ and $\mathcal{P}_2= \mathcal{Q}_2/\tilde{O}(\Lambda_h)$.
\begin{defi}
We denote by $\mathcal{M}_2$ the coarse moduli space of $\langle 2 \rangle$-polarized hyper-K\"ahler manifolds equivalent by deformation to the Hilbert square on a K3 surface.
\end{defi}
The moduli space $\mathcal{M}_2$ has dimension $20$ and is irreducible; it comes with a period map $\mathcal{M}_2\to \mathcal{P}_2$ which is an open embedding, see \cite[Theorem $1.10$]{markman2011}.

\medskip

There is an open embedding $\LG^0// \PGL(6) \hookrightarrow \mathcal{M}_2$, cf. \cite[Proposition $1.2.1$]{ogrady2011}, which induces a period map $\LG^0// \PGL(6)\to \mathcal{P}_2$.

We call $g=e-f$ the orthogonal complement of $h$ in the fixed copy of $U$ and $\ell$ a generator of the $\langle -2 \rangle$-part of $\Lambda$: we have $\Lambda_h=U^{\oplus 2}\oplus E_8(-1)^{\oplus 2}\oplus \IZ g\oplus \IZ\ell$.

\begin{defi}
We denote by $j\in O(\Lambda_h)$ the isometry exchanging $g$ and $\ell$. For simplicity, we also call $j$ the induced involution on $\mathcal{P}_2$.
\end{defi}

The discriminant of $\Lambda_h$ is isomorphic to $\IZ/2\IZ\times \IZ/2\IZ$ as a group. We identify the class of $g/2$ as $(1,0)$ and the class of $\ell/2$ as $(0,1)$. Since $\bar{j}((1,0))=(0,1)$, by \cite[Corollary $1.5.2$]{nikulin1979} the involution $j$ does not extend to an isometry of $\Lambda$.

Since the period map is an open embedding, the involution $j$ induces a birational involution on $\mathcal{M}_2$.

\begin{teo}\cite[Theorem $1.1$]{ogrady2006}
The involution induced on $\mathcal{M}_2$ by $j$ is the birational involution that sends $(X_A,D_A)$
to $(X_{A^\perp},D_{A^\perp})$. 
\end{teo}

The period map $(\LG-\Sigma )\to \mathcal{P}_2$ commutes with $j$ and 
the involution on $(\LG-\Sigma )//\PGL(6)$ sending $[A]$ to $[A^\perp]$, see \cite{ogrady2006a} and \cite{ogrady2010}.

\medskip

We consider now $A\in \Pi - (\Delta \cup \Sigma)$ and we fix $[V_5]\in Y^3_{A^\perp}$.
Let $S$ be the $GM$ variety associated to $(V_6,V_5,A)$; we focus on it as the $\langle 10\rangle$-polarized $K3$ surface $(S,H)$, with $H$ the polarization associated to the embedding in $\IP(W)$. 
The line bundle $H_2-2\delta$ is big and nef but not ample for a general choice of $A$.

The period map $\mathcal{M}_2\to \mathcal{P}_2$ extends to a surjective map $\mathcal{M}^\circ_2\to \mathcal{P}_2$, where $\mathcal{M}^\circ_2$ is the moduli space of hyper-K\"ahler manifolds with a big and nef divisor whose square is $2$.
Similarly, the period map on $\LG^0// \PGL(6)$
extends to a period map $(\LG - \Sigma)// \PGL(6)\to \mathcal{P}_2$.
The period point of $[A]$ is the image of the period point of $(S^{[2]}, H_2-2\delta)$ via the involution $j$, cf. \cite{ogrady2010} (but beware, $A$ is actually the dual of the Lagrangian subspace that O'Grady called $A$ in his paper).

As an easy consequence, it is possible to compute the Néron-Severi of $X_A$ from the one of $S$.

\begin{lemma}\label{pic di EPW e Hilb}
The involution $j$ induces an isometry $\NS(S^{[2]})\cap (H_2-2\delta)^\perp\cong \NS(X_A)\cap D_A^\perp$, in particular the Néron-Severi groups of $X_A$ and $S^{[2]}$ have same rank.
\end{lemma}

\bibliography{mabiblio}

\begin{bibdiv}
\begin{biblist}

\bib{debarre2014}{article}{
      author={Debarre, Olivier},
      author={Iliev, Atanas},
      author={Manivel, Laurent},
       title={Special prime {{Fano}} fourfolds of degree 10 and index 2},
        date={2014},
     journal={Recent advances in algebraic geometry},
      volume={417},
       pages={123},
}

\bib{debarre2015}{article}{
      author={Debarre, Olivier},
      author={Kuznetsov, Alexander},
       title={Gushel\textendash{{Mukai}} varieties: Classification and
  birationalities},
        date={2018},
     journal={Algebraic Geometry},
      volume={5},
       pages={15\ndash 76},
}

\bib{debarre2019a}{article}{
      author={Debarre, Olivier},
      author={Kuznetsov, Alexander},
       title={Double covers of quadratic degeneracy and {{Lagrangian}}
  intersection loci},
        date={2019},
     journal={Mathematische Annalen},
      volume={378},
       pages={1435–1469},
}

\bib{debarre2019}{article}{
      author={Debarre, Olivier},
      author={Kuznetsov, Alexander},
       title={Gushel\textendash{{Mukai}} varieties: Linear spaces and periods},
        date={2019},
     journal={Kyoto Journal of Mathematics},
      volume={59},
       pages={857–953},
}

\bib{debarre2020}{article}{
      author={Debarre, Olivier},
      author={Kuznetsov, Alexander},
       title={Gushel\textendash{{Mukai}} varieties: Intermediate
  {{Jacobians}}},
        date={2020},
     journal={\'Epijournal Géom. Algébrique},
      volume={4},
}

\bib{debarre2020a}{article}{
      author={Debarre, Olivier},
      author={Kuznetsov, Alexander},
       title={Gushel\textendash{{Mukai}} varieties: Moduli},
        date={2020},
     journal={International Journal of Mathematics},
       pages={2050013},
}

\bib{debarre2017}{article}{
      author={Debarre, Olivier},
      author={Macr{\`\i}, Emanuele},
       title={On the period map for polarized hyperk{\"a}hler fourfolds},
        date={2019},
     journal={International Mathematics Research Notices},
      volume={2019},
      number={22},
       pages={6887\ndash 6923},
}

\bib{debarre2021}{article}{
      author={Debarre, Olivier},
      author={Mongardi, Giovanni},
       title={Gushel--mukai varieties with many symmetries and an explicit
  irrational gushel--mukai threefold},
        date={2022},
     journal={Bollettino dell'Unione Matematica Italiana},
      volume={15},
       pages={133–161},
}

\bib{eisenbud2001}{article}{
      author={Eisenbud, David},
      author={Popescu, Sorin},
      author={Walter, Charles},
       title={Lagrangian subbundles and codimension 3 subcanonical subschemes},
        date={2001},
     journal={Duke Mathematical Journal},
      volume={107},
      number={3},
       pages={427\ndash 467},
}

\bib{ferretti2012}{article}{
      author={Ferretti, Andrea},
       title={The {{Chow}} ring of double {{EPW}} sextics},
        date={2012},
     journal={Algebra \& Number Theory},
      volume={6},
      number={3},
       pages={539\ndash 560},
}

\bib{garbagnati2013}{article}{
      author={Garbagnati, Alice},
       title={Elliptic k3 surfaces with abelian and dihedral groups of
  symplectic automorphisms},
        date={2013},
     journal={Communications in Algebra},
      volume={41},
      number={2},
       pages={583\ndash 616},
}

\bib{garbagnati2009}{article}{
      author={Garbagnati, Alice},
      author={Sarti, Alessandra},
       title={Elliptic fibrations and symplectic automorphisms on {{K3}}
  surfaces},
        date={2009},
     journal={Communications in Algebra},
      volume={37},
      number={10},
       pages={3601\ndash 3631},
}

\bib{greer2014}{article}{
      author={Greer, Francois},
      author={Li, Zhiyuan},
      author={Tian, Zhiyu},
       title={Picard groups on moduli of {{K3}} surfaces with {{Mukai}}
  models},
        date={2014},
     journal={International Mathematics Research Notices},
      volume={2015},
      number={16},
       pages={7238\ndash 7257},
}

\bib{gritsenko2010a}{article}{
      author={Gritsenko, Valery},
      author={Hulek, Klaus},
      author={Sankaran, Gregory~K.},
       title={Moduli spaces of irreducible symplectic manifolds},
        date={2010},
     journal={Compositio Mathematica},
      volume={146},
      number={2},
       pages={404\ndash 434},
}

\bib{gritsenko2010}{article}{
      author={Gritsenko, Valery},
      author={Hulek, Klaus},
      author={Sankaran, Gregory~K.},
       title={Moduli of {{K3}} surfaces and irreducible symplectic manifolds},
        date={2013},
     journal={Handbook of moduli. Vol. I, Adv. Lect. Math.},
      number={24},
       pages={459–526},
}

\bib{gritsenko2008}{article}{
      author={Gritsenko, Valery~A.},
      author={Hulek, Klaus},
      author={Sankaran, Gregory~K.},
       title={Abelianisation of orthogonal groups and the fundamental group of
  modular varieties},
        date={2009},
     journal={Journal of Algebra},
      volume={322},
      number={2},
       pages={463\ndash 478},
}

\bib{huybrechts2016}{book}{
      author={Huybrechts, Daniel},
       title={Lectures on {{K3}} surfaces},
   publisher={{Cambridge University Press}},
        date={2016},
      volume={158},
}

\bib{iliev2011}{inproceedings}{
      author={Iliev, Atanas},
      author={Manivel, Laurent},
       title={Fano manifolds of degree ten and {{EPW}} sextics},
        date={2011},
   booktitle={Annales scientifiques de l'{{\'Ecole Normale Sup\'erieure}}},
      volume={44},
       pages={393\ndash 426},
}

\bib{iskovskikh1999}{book}{
      author={Iskovskikh, V.~A.},
      author={Prokhorov, Yu~G.},
       title={Algebraic {{Geometry V}}. {{Fano Varieties}}, {{Encyclopedia}} of
  {{Math}}. {{Sci}}., vol. 47},
   publisher={{Springer}},
        date={1999},
}

\bib{knutsen2004}{book}{
      author={Knutsen, Andreas~L.},
      author={Johnsen, Trygve},
       title={K3 projective models in scrolls},
   publisher={{Springer}},
        date={2004},
}

\bib{markman2011}{incollection}{
      author={Markman, Eyal},
       title={A survey of {{Torelli}} and monodromy results for
  holomorphic-symplectic varieties},
        date={2011},
   booktitle={Complex and differential geometry},
   publisher={{Springer}},
       pages={257\ndash 322},
}

\bib{mecchia2011}{article}{
      author={Mecchia, Mattia},
      author={Zimmermann, Bruno},
       title={On finite groups acting on homology 4-spheres and finite
  subgroups of {{SO}} (5)},
        date={2011},
     journal={Topology and its Applications},
      volume={158},
      number={6},
       pages={741\ndash 747},
}

\bib{mongardi2013a}{article}{
      author={Mongardi, Giovanni},
       title={On symplectic automorphisms of hyper-kähler fourfolds of k3 [2]
  type},
        date={2013},
     journal={Michigan Mathematical Journal},
      volume={62},
       pages={537\ndash 550},
}

\bib{mongardi2016arx}{article}{
      author={Mongardi, Giovanni},
       title={Towards a classification of symplectic automorphisms on manifolds
  of {{K3}}[n]-type},
        date={2016},
     journal={Mathematische Zeitschrift},
      number={282},
       pages={651 \ndash  662},
}

\bib{mukai}{article}{
      author={Mukai, Shigeru},
       title={New developments in the theory of {{Fano}} threefolds: {{Vector}}
  bundle method and moduli problems, {{Sugaku}} 47 (1995), 124-144},
}

\bib{nikulin1979a}{article}{
      author={Nikulin, Viacheslav~Valentinovich},
       title={Finite groups of automorphisms of {{K\"ahlerian K}}\_3 surfaces},
        date={1979},
     journal={Trudy Moskovskogo Matematicheskogo Obshchestva},
      volume={38},
       pages={75\ndash 137},
}

\bib{nikulin1979}{article}{
      author={Nikulin, Viacheslav~Valentinovich},
       title={Integral symmetric bilinear forms and some of their
  applications},
        date={1979},
     journal={Izvestiya Rossiiskoi Akademii Nauk. Seriya Matematicheskaya},
      volume={43},
      number={1},
       pages={111\ndash 177},
}

\bib{ogrady2005}{article}{
      author={O'Grady, Kieran~G.},
       title={Involutions and linear systems on holomorphic symplectic
  manifolds},
        date={2005},
     journal={Geometric \& Functional Analysis GAFA},
      volume={15},
      number={6},
       pages={1223\ndash 1274},
}

\bib{ogrady2006a}{article}{
      author={O'Grady, Kieran~G.},
       title={Irreducible symplectic 4-folds and
  {{Eisenbud}}-{{Popescu}}-{{Walter}} sextics},
        date={2006},
     journal={Duke Mathematical Journal},
      volume={134},
      number={1},
       pages={99\ndash 137},
}

\bib{ogrady2006}{article}{
      author={O'Grady, Kieran~G.},
       title={Dual double {{EPW}}-sextics and their periods},
        date={2008},
     journal={Pure and Applied Mathematics Quarterly},
      volume={4},
      number={2},
       pages={427\ndash 468},
}

\bib{ogrady2012}{article}{
      author={O'Grady, Kieran~G.},
       title={{{EPW}}-sextics: Taxonomy},
        date={2012},
     journal={manuscripta mathematica},
      volume={138},
      number={1-2},
       pages={221\ndash 272},
}

\bib{ogrady2010}{article}{
      author={O'Grady, Kieran~G.},
       title={Double covers of {{EPW}}-sextics},
        date={2013},
     journal={Michigan Mathematical Journal},
      volume={62},
       pages={143–184},
}

\bib{ogrady2011}{article}{
      author={O’Grady, Kieran~G.},
       title={Moduli of double epw-sextics},
        date={2016},
      volume={240},
      number={1136},
}

\bib{piontkowski1999}{article}{
      author={Piontkowski, Jens},
      author={{Van de Ven}, Antonius},
       title={The automorphism group of linear sections of the {{Grassmannians
  G}} (1, {{N}})},
        date={1999},
     journal={Documenta Mathematica},
      volume={4},
       pages={623\ndash 664},
}

\bib{xiao1996a}{inproceedings}{
      author={Xiao, Gang},
       title={Galois covers between $ k3 $ surfaces},
        date={1996},
   booktitle={Annales de l'institut fourier},
      volume={46},
       pages={73\ndash 88},
}

\end{biblist}
\end{bibdiv}
\end{document}